\documentclass{amsart}
\usepackage{amssymb}\usepackage{amsmath}
\begin{document}
\newtheorem{theorem}{Theorem}[section]
\newtheorem{lemma}[theorem]{Lemma}
\newtheorem{remark}[theorem]{Remark}
\newtheorem{definition}[theorem]{Definition}
\newtheorem{corollary}[theorem]{Corollary}
\newtheorem{example}[theorem]{Example}
\def\Spec\{#1\}{\operatorname{Spec}\{#1\}}
\def\PSpec\{#1\}{\widetilde{\operatorname{Spec}}\{#1\}}
\def\proof{\medbreak\noindent{\it Proof.\ }}
\def\endproof{\hfill\qedbox\smallbreak\noindent}
\def\qedbox{\hbox{$\rlap{$\sqcap$}\sqcup$}}
\makeatletter
  \renewcommand{\theequation}{%
   \thesection.\alph{equation}}
  \@addtoreset{equation}{section}
 \makeatother
\def\MM{\mathfrak{M}}

\title[Projectively Osserman manifolds]
{Projectively Osserman manifolds}
\author[Brozos-V\'azquez et al.]
{M. Brozos-V\'azquez, P. Gilkey, S. Nik\v cevi\'c, and U. Simon}
\address{MBV: Department of Geometry and Topology, Faculty of
Mathematics, University of Santiago de Compostela, Santiago de Compostela
15782, Spain}
\email{mbrozos@usc.es}
\address{PG: Mathematics Department, University of Oregon, Eugene Or 97403,
USA}
\email{gilkey@uoregon.edu}
\address{SN: Mathematical Institute, Sanu, Knez Mihailova 35, p.p. 367, 11001 Belgrade, Serbia.}
\begin{email}{stanan@mi.sanu.ac.yu}\end{email}
\address
{US: Institut f\"ur Mathematik, Technische Universit\"at
Berlin\\ Strasse des 17. Juni 135, D-10623 Berlin, Deutschland}
\begin{email}{simon@math.tu-berlin.de}\end{email}

\begin{abstract}
One says that a smooth manifold $M$ of dimension $m$ is a pseudo-Riemannian manifold of signature $(p,q)$ if the tangent bundle $TM$
is equipped with a smooth non-degenerate symmetric inner product $g$ of signature $(p,q)$ where $p+q=m$. Similarly one says that $M$ is
an affine manifold if
$TM$ is equipped with a torsion free connection $\nabla$. One says $g$ is Osserman if the eigenvalues of the Jacobi operator are constant on the
pseudo-sphere bundles of unit timelike and spacelike vectors. We extend this
concept from the pseudo-Riemannian to the affine setting to define the notion of a projectively Osserman manifold. This notion is the focus
of the paper. We establish some basic results concerning projectively Osserman manifolds and exhibit examples of this structure
which arise in several different geometrical contexts.
\end{abstract}

\maketitle

\section{Introduction}\label{sect-1}

Let $\mathcal{M}=(M,g)$ be an $m$-dimensional pseudo-Riemannian manifold of
signature $(p,q)$ where $p+q=m$. Let $\nabla$ be the Levi-Civita connection defined by the metric $g$,
let $\mathcal{R}(x,y):=\nabla_x\nabla_y-\nabla_y\nabla_x-\nabla_{[x,y]}$ be the curvature operator, and let
$\mathcal{J}(x):y\rightarrow\mathcal{R}(y,x)x$ be the Jacobi operator. Let
$\Spec\{\mathcal{J}(x)\}\subset\mathbb{C}$ be the set of eigenvalues of $\mathcal{J}(x)$ and let
$S^\pm(\mathcal{M})$ be the pseudo-sphere bundles of unit spacelike (+) and unit timelike ($-$) tangent
vectors. One says that $\mathcal{M}$ is {\it spacelike Osserman at $P\in M$} if for every 
$x,y\in S^+(T_PM,g_P)$, $\Spec\{\mathcal{J}(x)\}=\Spec\{\mathcal{J}(y)\}$.
One says that $\mathcal{M}$ is pointwise spacelike Osserman if it is spacelike Osserman at every point of $M$
and that $\mathcal{M}$ is globally spacelike Osserman if the eigenvalue structure does not in fact depend on
the point in question. The notion {\it timelike Osserman} is defined by replacing
$S^+$ by
$S^-$ as appropriate. Note that if $p>0$ and if $q>0$, then work of Garc\'{\i}a--R\'{\i}o et al. 
\cite{GKVV99}  shows these are equivalent notions.
The investigation of Osserman manifolds has been an extremely active and fruitful one in recent years; we refer to
\cite{GKVm02,G01,G07} for further details.

In this paper, we wish to generalize these notions to the affine setting. Let $\mathcal{A}:=(M,\nabla)$
be an {\it affine manifold} where
$\nabla$ is a torsion free connection on $TM$. Again, let
$\mathcal{R}(x,y)$ be the curvature operator and let
$\mathcal{J}(x):y\rightarrow\mathcal{R}(y,x)x$ be the associated {\it Jacobi operator}; we will write $\mathcal{J}_\nabla$
when it is necessary to distinguish the role of the connection. Let
$\Spec\{\mathcal{J}(x)\}\subset\mathbb{C}$ be the spectrum of the Jacobi operator; since $\mathcal{J}(x)x=0$,
$0\in\operatorname{Spec}\{\mathcal{J}(x)\}$. Since
$\mathcal{J}(cx)=c^2\mathcal{J}(x)$,
$\Spec\{\mathcal{J}(cx)\}=c^2\Spec\{\mathcal{J}(x)\}$. In the pseudo-Riemannian
setting, we eliminated this rescaling effect by assuming that $g(x,x)=\pm1$. As this
normalization is not available in the affine setting, we must proceed slightly differently.

Recall that two non-zero points $u,v\in\mathbb{R}^m$ are said to be {\it projectively equivalent} if $u=cv$
for some $0\ne c\in\mathbb{R}$. This motivates the following definition; the role of $\{0\}$ is distinguished
and introduces a small amount of technical fuss.
\begin{definition}\label{defn-1}
\rm We say that an affine manifold $\mathcal{A}=(M,\nabla)$ is {\it projectively Osserman} at a point $P\in
M$ if there exists a subset $\mathcal{S}_P$ of $\mathbb{C}$ so that for any tangent vector $x\in
T_PM$, $\Spec\{\mathcal{J}(x)\}=c(x)\mathcal{S}_P$ for some suitably chosen complex number $c(x)$. We say
that
$\mathcal{A}$ is {\it pointwise projectively Osserman} if it is projectively Osserman at every point of $M$.
We say that
$\mathcal{A}$ is {\it globally projectively Osserman} if
$\mathcal{S}$ can be chosen independently of $P$.\end{definition}

This is related to earlier work by Garc\'{\i}a--R\'{\i}o et al.
\cite{GKVV99}. One says $\mathcal{A}=(M,\nabla)$ is {\it affine Osserman} if
$\Spec\{\mathcal{J}(x)\}=\{0\}$ for all tangent vectors $x$; such a manifold admits a natural  neutral signature Osserman metric, called the
Riemannian extension, on the cotangent bundle $T^*M$. Clearly any affine Osserman manifold
is projectively Osserman. The Riemannian extension is Osserman if and only if $\mathcal{A}$ is affine Osserman.

We have chosen to work with $\Spec\{\mathcal{J}(x)\}$, it is also possible to work with the unordered collection of eigenvalues
$\PSpec\{\mathcal{J}(x)\}$ where each eigenvalue is repeated according to multiplicity; working with $\PSpec\{\mathcal{J}(x)\}$ instead of
$\Spec\{\mathcal{J}(x)\}$ gives rise to the notion of {\it strongly projective Osserman}. Fortunately, these are equivalent concepts as we will
show in Section \ref{sect-1a}.

In this paper, we will exhibit several examples of projectively
Osserman affine manifolds. In Section
\ref{sect-2}, we discuss examples which arise from Osserman geometry. In Section \ref{sect-3}, we discuss projectively Osserman Walker
manifolds. In Section
\ref{sect-4}, we discuss examples from affine hypersurface theory. In Section \ref{sect-5}, we discuss the Weyl projective tensor.

\section{Eigenvalue multiplicities}\label{sect-1a}

This section is devoted to the proof of the following technical result.
\begin{lemma}\label{lem-2.1}
Let $\mathcal{A}=(M,\nabla)$ and let $P\in M$. Then $\mathcal{A}$ is projectively Osserman at $P$ if and only if $\mathcal{A}$
is strongly projectively Osserman at $P$.
\end{lemma}

\begin{proof} We must show that $\mathcal{A}$ is projectively Osserman at $P$ implies that $\mathcal{A}$ is strongly projectively Osserman at $P$
as the reverse implication is trivial. If
$\Spec\{\mathcal{J}(x)\}=\{0\}$ for all
$x\in T_PM$, then
$\PSpec\{\mathcal{J}(x)\}$ is the unordered set where $0$ is repeated with multiplicity $m$ and there is nothing to prove. We therefore
suppose there exists
$y\in T_PM$ so $\mathcal{J}(y)$ has a non-zero eigenvalue. As the eigenvalues vary continuously, as long as the spectrum does not degenerate to
$\{0\}$, the multiplicities are constant. Thus there is a small neighborhood
$\mathcal{O}$ of $y$ in $T_PM$ such that $\PSpec\{\mathcal{J}(x)\}=c(x)\PSpec\{\mathcal{J}(y)\}$ for $x\in\mathcal{O}$ where $c(x)\ne0$. Let
$p(x;t):=\det\{\mathcal{J}(x)-t\operatorname{id}\}$ be the characteristic polynomial.  Decompose 
$$p(x;t)=\prod_{i=1}^m(\lambda_i(x)-t)=\kappa_m(x)+\kappa_{m-1}(x)t+...+\kappa_kt^{m-k}+...+\kappa_0(x)t^m$$ where the
coefficients $\kappa_\nu(x)$ are the elementary symmetric functions of the eigenvalues $\{\lambda_1(x),...,\lambda_m(x)\}$. We note that
$\kappa_\nu$ is a polynomial of degree
$2k$ in the coordinate functions of
$x$ relative to some basis for $V$. Furthermore, if $x\in\mathcal{O}$, then $\lambda_\nu(x)=\kappa_\nu(x)\lambda_\nu(y)$ so 
\begin{equation}\label{eqn-2.a}
\kappa_\nu(x)=\kappa_\nu(y)c(x)^\nu\quad\text{for}\quad x\in\mathcal{O}\quad\text{and for}\quad 0\le\nu\le m\,.
\end{equation}
As $p(y;t)\ne(-t)^m$, there is $1\le\nu\le m$ so $\kappa_\nu(y)\ne0$. Thus
$$c(x):=\left\{\frac{\kappa_\nu(x)}{\kappa_\nu(y)}\right\}^{1/k}$$
is an analytic function of $x$. We may then complexify $V$ and $\mathcal{J}$ and consider the open dense subset $\mathcal{U}\subset
V\otimes\mathbb{C}$ where
$\Spec\{\mathcal{J}(x)\}\ne\{0\}$. We use analytic continuation to see that Equation (\ref{eqn-2.a}) holds for all $x\in\mathcal{U}\cap V$.
Consequently $p(y;t)=c(x)^{-m}p(x;c(x)t)$ if $\operatorname{Spec}\{\mathcal{J}(x)\}\ne\{0\}$ and the
desired result follows.
\end{proof}

\section{Osserman Manifolds}\label{sect-2}

One has the following observation:

\begin{theorem}\label{thm-2.1} Let $\mathcal{M}=(M,g)$ be a pseudo-Riemannian manifold. If
$\mathcal{M}$ is Osserman at $P\in M$, then $\mathcal{M}$ is projectively Osserman at $P$.
\end{theorem}

\begin{proof} Let $\mathcal{M}$ have signature $(p,q)$. This is immediate if the metric on $\mathcal{M}$ is positive definite since one has
$T_PM=\mathbb{R}\cdot S^+(T_PM,g_P)$ and since one also has $\mathcal{J}(cx)=c^2\mathcal{J}(x)$. The argument is the same
if the metric is negative definite and we therefore suppose $p>0$ and $q>0$. Let $\mathcal{M}$  be
spacelike Osserman at $P$ and let $\mathcal{S}_P^+:=\Spec\{\mathcal{J}(x)\}$ for any $x\in
S^+(T_PM,g_P)$. Let
$\varrho_\nu(x):=\operatorname{Tr}\{\mathcal{J}(x)^\nu\}$ for any
$x\in T_PM$. As
$\Spec\{\mathcal{J}(x)\}$ is constant on $S^+(T_PM,g_P)$, the eigenvalues and hence the
eigenvalue multiplicities are constant on $S^+(T_PM,g_P)$. This implies $\varrho_\nu(x)=\varrho_\nu$ is constant on
$S^+(T_PM,g_P)$. Since $\varrho_\nu(cx)=c^{2i}\varrho_\nu(x)$, we have that $\varrho_\nu(x)=g(x,x)^\nu\varrho_\nu$ if $x$ is spacelike. This
polynomial identity holds on an open subset of $T_PM$ and hence holds identically:
$$\operatorname{Tr}\{\mathcal{J}(x)^\nu\}=g(x,x)^\nu\varrho_\nu\quad\text{for all}\quad x\in T_PM\,.$$
It now follows, of course, that 
$\Spec\{\mathcal{J}(x)\}=g(x,x)\mathcal{S}_P^+$ for any $x\in T_PM$ and hence $\mathcal{M}$ is
projectively Osserman at $P$.
\end{proof}

Affine Osserman tensors play a central role:

\begin{theorem}\label{thm-2.2} Let $\mathcal{A}_1:=(M_1,\nabla_1)$ be projectively Osserman at $P_1\in M_1$
and let
$\mathcal{A}_2:=(M_2,\nabla_2)$ be affine Osserman at $P_2\in M_2$. Then the product structure
$\mathcal{A}:=(M_1\times M_2,\nabla_1\oplus\nabla_2)$ is projectively Osserman at
$P=(P_1,P_2)$.
\end{theorem}

\medbreak\noindent{\it Proof.} 
If $x=(x_1,x_2)\in T_{(P_1,P_2)}(M_1\times M_2)$, then $\mathcal{J}(x)=\mathcal{J}(x_1)\oplus\mathcal{J}(x_2)$
so
\begin{eqnarray*}
&&\Spec\{\mathcal{J}(x)\}
=\Spec\{\mathcal{J}(x_1)\}\cup\Spec\{\mathcal{J}(x_2)\}\\
&=&\Spec\{\mathcal{J}(x_1)\}\cup\{0\}=\Spec\{\mathcal{J}(x_1)\}
=c(x_1)\mathcal{S}_{P_1}\,.\qquad\qedbox
\end{eqnarray*}

We use this ansatz to construct new examples. Give the sphere $S^n$ and the torus $T^k$
the usual metrics where $n\ge2$ and $k\ge1$. Then $S^n$ is Osserman and $T^k$ is flat so $S^n\times T^k$ is
projectively Osserman. On the other hand, $S^n\times T^k$ is not Osserman. Thus there are projectively
Osserman manifolds which are not Osserman. Furthermore, while any flat manifold is affine Osserman, there are
other examples:
\begin{example}\label{ex-2.3}
\rm Follow the discussion in \cite{GIZ03}. Let $(x_1,...,x_p,\tilde x_1,...,\tilde x_p)$ be coordinates on
$\mathbb{R}^{2p}$. Let
$\psi$ be a smooth symmetric $2$-tensor field on
$\mathbb{R}^p$. Define a pseudo-Riemannian metric of neutral signature $(p,p)$ on $\mathbb{R}^{2p}$ whose
non-zero components are, up to the usual $\mathbb{Z}_2$-symmetries, given by:
$$g(\partial_{x_i},\partial_{x_j})=\psi_{ij}(x_1,...,x_p)\quad\text{and}\quad
  g(\partial_{x_i},\partial_{\bar x_j})=\delta_{ij}\,.$$
As $\mathcal{J}(x)^2=0$ for all $x$, $\Spec\{\mathcal{J}(x)\}=\{0\}$ for all $x$ as desired. 
\end{example}

\begin{example}\label{ex-2.4}
\rm Follow the discussion in \cite{GN04}. For $s\ge2$, choose coordinates
$$(u_1,...,u_s,t_1,...,t_s,v_1,...,v_s)$$ 
 on $\mathbb{R}^{3s}$. Let $f_i\in
C^\infty(\mathbb{R})$ be given. Define a
pseudo-Riemannian metric
$g$ of signature
$(2s,s)$ on $\mathbb{R}^{3s}$ whose non-zero components are given by:
\begin{eqnarray*}
g(\partial_{u_i},\partial_{u_i})=-2(f_1(u_1)+u_1t_1+...+f_s(u_s)+u_st_s),\\
 g(\partial_{u_i},\partial_{v_i})=g(\partial_{v_i},\partial_{u_i})=1,\quad\text{and}\quad
 g(\partial_{t_i},\partial_{t_i})=-1\,.
\end{eqnarray*}
Then $\mathcal{J}(x)^3=0$ for all $x$ so $\mathcal{A}:=(\mathbb{R}^{3s},\nabla)$ is affine
Osserman.
\end{example}

One says that an affine manifold $\mathcal{A}=(M,\nabla)$ is $k$-affine curvature homogeneous if given any two points $P,Q\in M$, there
is an isomorphism $\phi:T_PM\rightarrow T_QM$ so that $\phi^*\{\nabla^i\mathcal{R}_P\}=\nabla^i\mathcal{R}_Q$ for $0\le i\le k$. 

\begin{example}\label{ex-2.5}
\rm Follow the discussion in \cite{GN04d}.
Let $(x,y,z_0,...,z_\ell,\tilde x,\tilde y,\tilde z_0,...,\tilde z_\ell)$ be coordinates on $\mathbb{R}^{6+2\ell}$. Let $g_f$ be the
pseudo-Riemannian neutral signature metric defined by:
\begin{eqnarray*}
&&g_f(\partial_x,\partial_{\tilde x})=g_f(\partial_y,\partial_{\tilde y})=g_f(\partial_{z_i},\partial_{\tilde z_i})=1,\\
&&g_f(\partial_x,\partial_x)=f(y)+yz_0+...+y^{\ell+1}z_\ell\,.
\end{eqnarray*}
Assume $f^{(\ell+3)}>0$ and $f^{(\ell+4)}>0$. Let $\mathcal{A}_f:=(\mathbb{R}^{6+2\ell},\nabla_f)$ where $\nabla_f$ is the
Levi-Civita connection defined by $g_f$. Then $\mathcal{A}_f$ is affine Osserman and $\mathcal{A}_f$ is $(\ell+2)$-affine curvature homogeneous.
Furthermore, 
$\mathcal{A}_f$ is $(\ell+3)$-affine curvature homogeneous if and only if $f^{(\ell+3)}(y)=ae^{by}$ for $a>0$ and $b>0$ real constants; this
happens if and only if
$\mathcal{A}_f$ is affine homogeneous.
\end{example}

The manifolds in Examples \ref{ex-2.3}, \ref{ex-2.4}, and \ref{ex-2.5} are generalized plane wave manifolds and hence are geodesically
complete. We refer to \cite{G07} for other examples and to \cite{GN05} for a further discussion of generalized plane wave manifolds.

\section{Walker manifolds}\label{sect-3}
The following family will be crucial for our study.

\begin{definition}\label{defn-3.1}
\rm Let $(x_1,x_2,x_3,x_4)$ be coordinates on $\mathbb{R}^4$. Consider the
following {\it Walker manifold} $\mathcal{M}:=(\mathbb{R}^4,g)$ of signature $(2,2)$ where
$$g(\partial_{x_1},\partial_{x_3})=g(\partial_{x_2},\partial_{x_4})=1\quad\text{and}\quad
  g(\partial_{x_3},\partial_{x_4})=g_{34}(x_1,x_2,x_3,x_4)\,.$$
\end{definition}

Let $\mathcal{W}$ be the {\it Weyl conformal curvature operator} of a pseudo-Riemannian manifold and let
$\mathcal{J}_W(x):y\rightarrow\mathcal{W}(y,x)x$ be the {\it conformal Jacobi operator}.  One says that $\mathcal{M}$ is
{\it pointwise conformally Osserman} if $\mathcal{J}_W$ has constant eigenvalues on $S^\pm(T_PM,g_P)$ for every point $P\in M$.
We showed \cite{BBGS05} that this is a conformal notion; $(M,g)$ is pointwise conformally Osserman if and only if $(M,e^{h}g)$
is pointwise conformally Osserman for any $h\in C^\infty(M)$.

Recall that in dimension
$4$,
$\mathcal{M}$ is conformally Osserman if and only if
$\mathcal{M}$ is either self-dual or anti-self-dual. We take the orientation
$dx_1dx_2dx_3dx_4$ for $\mathbb{R}^4$. If $f=f(x_1,x_2,x_3,x_4)$, let $f_{/i}:=\partial_{x_i}f$ and
$f_{/ij}:\partial_{x_i}\partial_{x_j}f$. One has the following result concerning these
manifolds
\cite{BGVa05, BGGV06a}:

\begin{theorem}\label{thm-3.2}
Let $\mathcal{M}$ be as in Definition \ref{defn-3.1}. Then
\begin{enumerate}
\item $\mathcal{M}$ is self-dual if and only if
$g_{34}=x_1p(x_3,x_4)+x_2q(x_3,x_4)+s(x_3,x_4)$.
\item $\mathcal{M}$ is anti-self-dual if
and only if
$g_{34}=x_1p(x_3,x_4)+x_2q(x_3,x_4)+s(x_3,x_4)$ $+\xi(x_1,x_4)+\eta(x_2,x_3)$ for
$p_{/3}=q_{/4}$ and $g_{34}p_{/3}-x_1p_{/34}-x_2p_{/33}-s_{/34}=0$.
\item The following assertions are equivalent:
\begin{enumerate}
\item $\mathcal{M}$ is Osserman.
\item $\mathcal{M}$ is Einstein.
\item The Ricci tensor is zero.
\item $g_{34}=x_1p(x_3,x_4)+x_2q(x_3,x_4)+s(x_3,x_4)$ where $p^2=2p_{/4}$, $q^2=2q_{/3}$, and
$pq=p_{/3}+q_{/4}$.
\end{enumerate}\end{enumerate}
\end{theorem}

\begin{remark}\label{rmk-3.3}
\rm Let $\mathcal{M}$ be as in Definition \ref{defn-3.1}. Results of \cite{BGGV06a} show that if $\mathcal{M}$ is Einstein, then
$g_{34}=x_1p(x_3,x_4)+x_2q(x_3,x_4)+s(x_3,x_4)$ where $p$ and $q$ have one of the following forms:
\begin{enumerate}
\item $p=q=0$.
\item $p=0$ and $q=-2(x_3+b(x_4))^{-1}$. 
\item $p=-2(x_4+a(x_3))^{-1}$ and $q=0$.
\item $p=-2(x_4+a)^{-1}$ and $q=-2(x_3+b)^{-1}$.
\item $p=-2(x_3+b_0+b_1x_4)^{-1}$ and $q=-2(x_4+a_0+a_1x_3)^{-1}$ where $a_1b_1=1$ and $a_0=b_0a_1$.
\end{enumerate}\end{remark}

 We now come to the main result of this section:
\begin{theorem}\label{thm-3.4}
Let $\mathcal{M}$ be as in Definition \ref{defn-3.1}. Then
$\mathcal{M}$
is globally projectively Osserman if and only if at least one of the following conditions holds:
\begin{enumerate}
\item $g_{34}=p(x_1,x_4)+s(x_3)$.
\item $g_{34}=q(x_2,x_3)+s(x_4)$.
\item $g_{34}=x_1p(x_3,x_4)+x_2q(x_3,x_4)+s(x_3,x_4)$.
\end{enumerate}
\end{theorem} 

\begin{proof} This is a computer assisted computation.
We begin by verifying the defining functions of Theorem \ref{thm-3.4} (1)-(3) define projectively
Osserman manifolds. We set $\xi:=\sum_iv_i\partial_{x_i}$. Suppose that
$$g_{34}=x_1p(x_3,x_4)+x_2q(x_3,x_4)+s(x_3,x_4)\,.$$
Then $\Spec\{\mathcal{J}(\xi)\}=\{0,a(\xi)\}$ where $0$ and $a(\xi)$ appear with multiplicity $2$ with
\begin{eqnarray*}
a(\xi)&=&\textstyle\frac14\{-v_4^2p(x_3,x_4)^2+2
v_3
v_4p(x_3,x_4) q(x_3,x_4) -v_3^2q(x_3,x_4) ^2\\
&&\qquad+2 v_4^2p_{/4}(x_3,x_4) -2 v_3
v_4q_{/4}(x_3,x_4) -2 v_3 v_4p_{/3}(x_3,x_4) \\
&&\qquad+2 v_3^2q_{/3}(x_3,x_4)\}\,.
\end{eqnarray*}

Next suppose that $g_{34}=p(x_1,x_4)+s(x_3)$; the case $g_{34}=q(x_2,x_3)+s(x_4)$ being analogous. One has
$\Spec\{\mathcal{J}(\xi)\}=\{0,a(\xi)\}$ where $0$ and $a(\xi)$ appear with multiplicity $2$ with
\begin{eqnarray*}
a(\xi)&=&-\textstyle\frac{1}{4} v_4 \{v_4 p_{/1}(x_1,x_4) ^2
-2 v_4 p_{/14}(x_1,x_4) -2 {v_1} p_{/11}(x_1,x_4)\}\,.
\end{eqnarray*}

Conversely, suppose $\mathcal{M}$ is as given in Definition \ref{defn-3.1} and that $\mathcal{M}$ is projectively Osserman. 
We suppose first that $0$ is at least a double eigenvalue. The coefficients of $\lambda v_1^3v_3^2v_4$ and of $\lambda v_1v_2^2v_3^2v_4$ in
the characteristic polynomial are seen to be
$$\textstyle\frac12g_{34/11}(g_{34/12}^2-g_{34/11}g_{34/22})
\quad\text{and}\quad \textstyle\frac12g_{34/22}(g_{34/12}^2-g_{34/11}g_{34/22})\,.$$
To ensure $0$ is at least a double eigenvalue, we set these two terms to zero and obtain two cases:
\smallbreak\noindent{\bf Case 1:} $(g_{34/12}^2-g_{34/11}g_{34/22})\ne0$. We then have $g_{34/11}=g_{34/22}=0$ so
$$g_{34}=x_1p(x_3,x_4)+x_2q(x_3,x_4)+x_1x_2r(x_3,x_4)+s(x_3,x_4)\,.$$
One computes the coefficient of $\lambda v_1v_2v_3^3v_4x_1^2$ in the characteristic polynomial to be $-\frac12r(x_3,x_4)^4$. Setting this
to zero implies $r=0$ so
$g_{34}$ is affine in $\{x_1,x_2\}$ which is one of the possibilities enumerated in the Theorem.

\medbreak\noindent{\bf Case 2:} We have the identity $g_{34/12}^2=g_{34/11}g_{34/22}$. Only the first and second derivatives appear in
the calculation of the curvature tensor. Thus we may approximate by the second order Taylor polynomial. For simplicity, we suppose the
point in question to be $x_1=x_2=0$. Set
\begin{eqnarray*}
g_{34}&=&a_0(x_3,x_4)+x_1a_1(x_3,x_4)+x_2a_2(x_3,x_4)\\
&+&x_1^2a_{11}(x_3,x_4)^2+x_2^2a_{22}(x_3,x_4)^2
+4x_1x_2a_{11}(x_3,x_4)a_{22}(x_3,x_4)\,.
\end{eqnarray*}
When $x_1=x_2=0$, the coefficient of $\lambda v_1v_2^2v_3^2v_4$ in the characteristic polynomial is seen to be
$12a_{11}(x_3,x_4)^2a_{22}(x_3,x_4)^4$. Thus we have $g_{34/12}=0$ and, without loss of generality $g_{34/11}=0$.
This means that
$$g_{34}=x_1p(x_3,x_4)+q(x_2,x_3,x_4)\,.$$
Setting the coefficient of $\lambda$ in the characteristic polynomial to zero then leads to:
\begin{eqnarray*}
0&=&\textstyle-\frac{1}{8}\{-2 p_{/3}(x_3,x_4)+p(x_3,x_4)q_{/2}(x_2,x_3,x_4)\}^2 q_{/22}(x_2,x_3,x_4),\\
0&=&\textstyle\frac{1}{4}\{p(x_3,x_4)^2-2 p_{/4}(x_3,x_4)\}\{-2p_{/3}(x_3,x_4)+p(x_3,x_4)q_{/2}(x_2,x_3,x_4)\}\\
&&\times q_{/22}(x_2,x_3,x_4),\\
0&=&\textstyle-\frac{1}{8}\{p(x_3,x_4)^2-2 p_{/4}(x_3,x_4)\}^2q_{/22}(x_2,x_3,x_4)\,.
\end{eqnarray*}
One possibility is $q_{/22}=0$. This implies $g_{34}$ is affine in $\{x_1,x_2\}$ as desired. The other possibility is
$q_{/22}\ne0$ so
$2p_{/3}(x_3,x_4)=p(x_3,x_4)q_{/2}(x_2,x_3,x_4)$.
Differentiating this relation with respect to $\partial_{x_2}$  yields $0=pq_{/22}$ so $p=0$ and
$$g_{34}=q(x_2,x_3,x_4)\,.$$
Zero is at least a double eigenvalue and the other
eigenvalues are
\begin{eqnarray*}
\lambda_{\pm}(x_2,x_3,x_4)&=&\textstyle\frac{1}{4} (-v_3^2 q_{/2}^2-2 v_3 v_4
q_{/24}+2
v_3^2 q_{/23}+
2 {v_2} v_3 q_{/22})\\&\pm& \sqrt{ -v_3^3 v_4 q_{/34}q_{/22}-{v_2}
v_3^2 v_4 q_{/24} q_{/22}}\,.\end{eqnarray*}

Since $q_{/22}\ne0$, $q_{/2}\ne0$. Setting $v_4=0$ and $v_3=2$ yields $\lambda_\pm=-q_{/2}^2\ne0$. By Lemma \ref{lem-2.1},
the remaining eigenvalue must be always be double eigenvalue so the square root
must vanish identically. Consequently
$q_{/34}=q_{/24}=0$ and thus
$g_{34}=q(x_2,x_3)+s(x_4)$ which is the possibility given in (2) of the Theorem.

We complete the proof by analyzing what happens if there is a tangent vector $\xi_0$ so that $0$ is a simple eigenvalue of $\mathcal{J}(\xi_0)$.
Lemma \ref{lem-2.1} then shows that for any tangent vector $\xi$ that either
$\operatorname{Spec}\{\mathcal{J}(\xi)\}=\{0\}$ or $0$ is a simple eigenvalue of $\mathcal{J}(\xi)$. We specialize and set
$v_4=0$. Zero is then a double eigenvalue and the other eigenvalue $a(\xi)$ is double eigenvalue as well where:
\begin{eqnarray*}
a(\xi)&=&\textstyle\frac{1}{4} v_3 \{-v_3 g_{34/2}(x_1,x_2,x_3,x_4) ^2+2 v_3
g_{34/23}(x_1,x_2,x_3,x_4) \\
&&\quad+2{v_2}
g_{34/22}(x_1,x_2,x_3,x_4) +2{v_1} g_{34/12}(x_1,x_2,x_3,x_4)\}
\,.
\end{eqnarray*}
Consequently, $a(\xi)=0$. 
This implies $g_{34/12}=g_{34/22}=0$; specializing to set $v_3=0$ instead of $v_4=0$ yields $g_{34/11}=0$ as well. Thus $g_{34}$ is
affine in
$(x_1,x_2)$.
\end{proof}

\section{Relative hypersurfaces}\label{sect-4}
We refer to   \cite{LSZ93, SVV91} for further material concerning the theory of relative hypersurfaces. 
Let $\mathcal{A}=(M,\nabla)$.
We begin with a technical observation. Suppose there is a quadratic form $\tau$ so that
\begin{equation}\label{eqn-4.a}
\mathcal{J}_\nabla(x)y=\tau(x,x)y-\tau(x,y)x\,.
\end{equation}
Because $\mathcal{J}_\nabla(x)x=0$, $0\in\Spec\{\mathcal{J}_\nabla(x)\}$. Suppose that $\tau(x,x)=0$.
One then has that $\mathcal{J}_\nabla(x)y=-\tau(x,y)x$ and thus $\mathcal{J}_\nabla(x)^2y=0$. Since $\mathcal{J}_\nabla(x)$ is nilpotent,
$\Spec\{\mathcal{J}_\nabla(x)\}=\{0\}$. On the other hand, suppose $\tau(x,x)\ne0$. Let $y\in x^\perp$, i.e. $\tau(x,y)=0$. Then
$\mathcal{J}(x)y=\tau(x,x)y$ and we conclude $\Spec\{\mathcal{J}_\nabla(x)\}=\{0,\tau(x,x)\}$. Thus $\nabla$ is
projectively Osserman by Equation (\ref{eqn-4.a}) as
$\operatorname{Spec}\{\mathcal{J}_\nabla\}=\{0,\tau(x,x)\}$.

Fix a relative normalization for a hypersurface $M$ in affine space $\mathbb{A}^{m+1}$. Let $(\nabla,g,\nabla^*)$
be the induced connection, the relative metric, and the conormal connection; this forms a {\it conjugate triple} and one has that:
$$xg(y,z)=g(\nabla_xy,z)+g(y,\nabla^*_xz)\,.$$
Let $\mathcal{R}$ and $\mathcal{R}^*$ be
the associated curvature operators, let $\operatorname{Ric}^*$ be the Ricci tensor of $\nabla^*$, and let $\rho^*$ be the Ricci operator of
$\nabla^*$. One then has
\begin{eqnarray*}
&&\mathcal{R}^*(v_1,v_2)w=\textstyle\frac1{m-1}\{\operatorname{Ric}^*(v_2,w)v_1-\operatorname{Ric}^*(v_1,w)v_2\},\\
&&\mathcal{R}(v_1,v_2)w=\textstyle\frac1{m-1}\{g(v_2,w)\rho^*v_1-g(v_1,w)\rho^*v_2\},\\
&&\mathcal{J}_{\nabla^*}y=\mathcal{R}^*(y,x)x=\textstyle\frac1{m-1}\left\{\operatorname{Ric}^*(x,x)y-\operatorname{Ric}^*(y,x)x\right\},\\
&&\mathcal{J}_{\nabla}y=\mathcal{R}(y,x)x=\textstyle\frac1{m-1}\left\{g(x,x)\rho^*y-g(y,x)\rho^*x\right\}\,.
\end{eqnarray*}

Let $H_P$ be the relative mean curvature; $M$ is said to be a {\it relative umbilic} at a point $P$ of $M$ if
$\operatorname{Ric}^*=(m-1)H_Pg$ at $P$. In
 this setting, $\operatorname{Ric}^*= \operatorname{Ric}$ so
$$\mathcal{J}_{\nabla^*}(x)y=\mathcal{J}_\nabla(x)y=H_P\left\{g(x,x)y-g(x,y)x\right\}$$
and consequently both $(M,\nabla)$ and $(M,\nabla^*)$ are  projectively Osserman at  $P$. 
If all points on $M$ are relative umbilics, then the hypersurface with
its relative normalization is said to be a {\it relative sphere}; $(M,\nabla)$ and $(M,\nabla^*)$ are  globally projectively Osserman in such a
situation. In particular: the spectra of the two connections, which are in general different, coincide on $M$.

\subsection{Centroaffine geometry} The centroaffine normalization of an affine hypersurface is given by a 
transversal position vector; this  is a relative normalization, and in terms of this relative normalization,
any centroaffine hypersurface is a relative sphere. Thus, for all centroaffine
hypersurfaces, the discussion of the last section applies.

\subsection{Examples} Affine spheres in the unimodular hypersurface theory,
relative spheres 
and centroaffine hypersurfaces are very big classes of hypersurfaces
that play an important role in the theory. Thus the foregoing sections give
further examples having the projectively Osserman property either pointwise or globally.

\section{The Weyl projective tensor}\label{sect-5}
There is no analogue of the Weyl conformal tensor in the geometry of
a manifold $M$ equipped with an affine connection $\nabla$, and thus we can not speak of conformally Osserman in this context. Instead, a
similar role is played by the Weyl projective curvature tensor which  is an invariant of the projective class generated by $\nabla$.

Recall that two connections
$\nabla$ and
$\nabla^{\sharp}$ are said to be {\it projectively equivalent} if there exists a $1$-form $\theta$ so that
\begin{equation}\label{eqn-5.a}
\nabla_uv-\nabla^{\sharp}_uv=\theta(u)v+\theta(v)u\,.
\end{equation}
We remark that the unparametrized geodesics of $\nabla$ and $\nabla^{\sharp}$ coincide if and only if both connections  are projectively equivalent. Thus the
algebraic definition given in Equation (\ref{eqn-5.a}) has great geometric significance.

The {\it projective curvature operator}
$\mathcal{P}=\mathcal{P}_\nabla$ is defined by:
\begin{eqnarray*}
 &&\mathcal{P}(v_1,v_2)w:=\mathcal{R}(v_1,v_2)w+\textstyle\frac1{m^2-1}\left\{m\operatorname{Ric}(v_1,w)v_2
 + \operatorname{Ric}(w,v_1)v_2\right\}\\
 &&\quad-\textstyle\frac1{m^2-1}\left\{m\operatorname{Ric}(v_2,w)v_1 
 +\operatorname{Ric}(w,v_2)v_1\right\}
 \\&&\quad+\textstyle
 \frac1{m+1}\left\{\operatorname{Ric}(v_1,v_2)w-\operatorname{Ric}(v_2,v_1)w\right\}\,.
\end{eqnarray*}
If $\nabla$ is Ricci symmetric, we have a somewhat simpler form:
\begin{eqnarray*}
\mathcal{P}(v_1,v_2)w:=\mathcal{R}(v_1,v_2)w-\textstyle\frac1{m-1}\left\{\operatorname{Ric}(v_2,w)v_1-\operatorname{Ric}(v_1,w)v_2\right\}\,.
\end{eqnarray*}
This operator satisfies the identities:
\begin{eqnarray*}
&&\mathcal{P}(v_1,v_2)=-\mathcal{P}(v_2,v_1),\quad\text{and}\\
&&\mathcal{P}(v_1,v_2)v_3+\mathcal{P}(v_2,v_3)v_1+\mathcal{P}(v_3,v_1)v_2=0\,.
\end{eqnarray*}
If $\nabla$ and $\nabla^\sharp$ are projectively equivalent and Ricci symmetric, then $\mathcal{P}_\nabla=\mathcal{P}_{\nabla^\sharp}$.

We use $\mathcal{P}$ to define the {\it projective Jacobi operator} $\mathcal{J}_P(x):y\rightarrow\mathcal{P}(y,x)x$; the fundamental observation
is then that if $\nabla$ is projectively equivalent to $\nabla^{\sharp}$, then 
$$\mathcal{J}_{P,\nabla}(x)=\mathcal{J}_{P,\nabla^{\sharp}}(x)\quad\text{for all}\quad x\,.$$

We note that we are using the word {\it projective} in two different settings. We shall say that $(M,\nabla)$ is
{\it pointwise projectively Weyl Osserman} if for every point $P\in M$, there is a subset $\mathcal{S}_P\subset\mathbb{C}$ so that
for any tangent vector $x\in T_PM$, $\Spec\{\mathcal{J}_{P,\nabla}(x)\}=c(x)\mathcal{S}_P$ for a suitably chosen constant $c(x)$. This notion
plays the role in the affine setting that the notion `conformal Osserman' plays in the geometric setting. In particular, any projectively flat
manifold is necessarily projectively Weyl Osserman.

\section*{Acknowledgments}
Research of M. Brozos-V\'azquez was partially supported by the Max
Planck Institute for Mathematics in the Sciences (Germany), by FPU
grant and by projects PGIDIT06PXIB207054PR and MTM2006-01432
(Spain). Research of P. Gilkey was partially supported by the Max Planck Institute in the Mathematical
Sciences (Leipzig, Germany), and by Project MTM2006-01432 (Spain). Research of S. Nik\v cevi\'c was partially
supported by the TU Berlin, and by Project 144032 (Serbia).  Research of 
 U. Simon was partially supported by DFG PI 158/4-5.

  \end{document}